\documentclass[11pt]{article}

\usepackage{amssymb,amsbsy,amsthm,amsmath,amscd,epsfig,graphicx,times}

\begin{document}

\title{\textbf{Deducing the Density Hales-Jewett Theorem from an infinitary removal lemma}}
\author{Tim Austin}
\date{}

\maketitle


\newenvironment{nmath}{\begin{center}\begin{math}}{\end{math}\end{center}}

\newtheorem{thm}{Theorem}[section]
\newtheorem*{thm*}{Theorem}
\newtheorem{lem}[thm]{Lemma}
\newtheorem{prop}[thm]{Propoisition}
\newtheorem{cor}[thm]{Corollary}
\newtheorem*{conj*}{Conjecture}
\newtheorem{dfn}[thm]{Definition}
\newtheorem{ques}[thm]{Question}
\theoremstyle{remark}


\newcommand{\A}{\mathcal{A}}
\newcommand{\B}{\mathcal{B}}
\newcommand{\I}{\mathcal{I}}
\newcommand{\frH}{\mathfrak{H}}
\renewcommand{\Pr}{\mathrm{Pr}}
\newcommand{\s}{\sigma}
\renewcommand{\P}{\mathcal{P}}
\renewcommand{\O}{\Omega}
\renewcommand{\o}{\omega}
\renewcommand{\S}{\Sigma}
\newcommand{\T}{\mathrm{T}}
\newcommand{\co}{\mathrm{co}}
\newcommand{\e}{\mathrm{e}}
\newcommand{\eps}{\varepsilon}
\renewcommand{\d}{\mathrm{d}}
\newcommand{\im}{\mathrm{i}}
\renewcommand{\l}{\lambda}
\newcommand{\U}{\mathcal{U}}
\newcommand{\G}{\Gamma}
\newcommand{\g}{\gamma}
\newcommand{\calL}{\mathcal{L}}
\renewcommand{\L}{\Lambda}
\newcommand{\hcf}{\mathrm{hcf}}
\newcommand{\FLat}{\mathrm{FLat}}
\newcommand{\F}{\mathcal{F}}
\renewcommand{\a}{\alpha}
\newcommand{\bbN}{\mathbb{N}}
\newcommand{\bbR}{\mathbb{R}}
\newcommand{\bbZ}{\mathbb{Z}}
\newcommand{\bbQ}{\mathbb{Q}}
\newcommand{\bbT}{\mathbb{T}}
\newcommand{\sfE}{\mathsf{E}}
\newcommand{\sfP}{\mathsf{P}}
\newcommand{\id}{\mathrm{id}}
\newcommand{\bb}[1]{\mathbb{#1}}
\newcommand{\fr}[1]{\mathfrak{#1}}
\renewcommand{\bf}[1]{\mathbf{#1}}
\renewcommand{\rm}[1]{\mathrm{#1}}
\renewcommand{\cal}[1]{\mathcal{#1}}
\newcommand{\fin}{\nolinebreak\hspace{\stretch{1}}$\lhd$}
\newcommand{\uhr}{\!\!\upharpoonright}

\newcommand{\into}{\hookrightarrow}

\renewcommand{\line}{\mathrm{line}}

\parskip 7pt

\parindent 0pt

\begin{abstract}
We offer a new proof of Furstenberg and Katznelson's density version
of the Hales-Jewett Theorem:
\begin{thm*}
For any $\delta > 0$ there is some $N_0 \geq 1$ such that whenever
$A \subseteq [k]^N$ with $N \geq N_0$ and $|A|\geq \delta k^N$, $A$
contains a \textbf{combinatorial line}: that is, for some $I
\subseteq [N]$ nonempty and $w_0 \in [k]^{[N]\setminus I}$ we have
\[A \supseteq \{w:\ w|_{[N]\setminus I} = w_0,\,w|_I = \rm{const.}\}.\]
\end{thm*}

Following Furstenberg and Katznelson, we first show that this result
is equivalent to a `multiple recurrence' assertion for a class of
probability measures enjoying a certain kind of stationarity.
However, we then give a quite different proof of this latter
assertion through a reduction to an infinitary removal lemma in the
spirit of Tao~\cite{Tao07} (and also its recent re-interpretation
in~\cite{Aus--newmultiSzem}).  This reduction is based on a
structural analysis of these stationary laws closely analogous to
the classical representation theorems for various partial
exchangeable stochastic processes in the sense of
Hoover~\cite{Hoo79}, Aldous~\cite{Ald82,Ald85} and
Kallenberg~\cite{Kal92}. However, the underlying combinatorial
arguments used to prove this theorem are rather different from those
required to work with exchangeable arrays, and involve crucially an
observation that arose during ongoing work by a collaborative team
of authors~\cite{Gow(online)} to give a purely finitary proof of the
above theorem.
\end{abstract}

\tableofcontents

\section{Introduction}

In this note we record a new proof of the following result of
Furstenberg and Katznelson:

\begin{thm}[Density Hales-Jewett Theorem]\label{thm:DHJ}
For any $\delta > 0$ there is some $N_0 \geq 1$ such that whenever
$A \subseteq [k]^N$ with $N \geq N_0$ and $|A|\geq \delta k^N$, $A$
contains a \textbf{combinatorial line}: that is, for some $I
\subseteq [N]$ nonempty and $w_0 \in [k]^{[N]\setminus I}$ we have
\[A \supseteq \{w:\ w|_{[N]\setminus I} = w_0,\,w|_I = \rm{const.}\}.\]
\end{thm}

This is the `density' version of the classical Hales-Jewett Theorem
of colouring Ramsey Theory.  The Hales-Jewett Theorem is one of the
central results of Ramsey Theory, partly because many other results
in that area can be deduced from it (see, for example, Chapter 2 of
Graham, Rothschild and Spencer~\cite{GraRotSpe90}).  Likewise, the
above density variant generalizes many other results in density
Ramsey Theory, such as the famous theorem of
Szemer\'edi~(\cite{Sze75}) and its multidimensional analog, also
proved by Furstenberg and Katznelson~(\cite{FurKat78}).

Following Furstenberg's discovery in~\cite{Fur77} of an alternative
proof of Szemer\'edi's Theorem via a conversion to a result in
ergodic theory, the use of ergodic-theoretic methods to prove
results in density Ramsey Theory has become widespread and powerful
(see, for example, the survey of Bergelson~\cite{Ber96}), and the
above result was one of the furthest-reaching consequences of this
program. In this paper we follow Furstenberg and Katznelson as far
as their re-interpretation of the above result in terms of
stochastic processes, but then we offer a new proof of that version
of the result.

In order to state the result about stochastic processes into which
Furstenberg and Katznelson convert Theorem~\ref{thm:DHJ}, let us
first define $[k]^\o$ to be the \textbf{infinite-dimensional
combinatorial space} over the alphabet $k$:
\[[k]^\o := \bigcup_{n\geq 1}[k]^n.\]

\begin{thm}[Infinitary Density Hales-Jewett Theorem]\label{thm:inf-DHJ}
For any $\delta > 0$, if $\mu$ is a Borel probability measure on
$\{0,1\}^{[k]^\o}$ for which
\[\mu\{\bf{x} \in \{0,1\}^{[k]^\o}:\ x_w = 1\} \geq \delta\quad\quad\forall w \in \{0,1\}^{[k]^\o}\]
then there are some $N\geq 1$, $I\subseteq [N]$ nonempty and $w_0
\in [k]^{[N]\setminus I}$ such that
\[\mu\{\bf{x} \in \{0,1\}^{[k]^\o}:\ x_w = 1\ \forall w \in [k]^N\ \rm{s.t.}\ w|_{[N]\setminus I} = w_0,\,w|_I = \rm{const.}\} > 0.\]
\end{thm}

It is shown in Proposition 2.1 of~\cite{FurKat91} that
Theorems~\ref{thm:DHJ} and~\ref{thm:inf-DHJ} are equivalent.  Here
we will assume this first step of Furstenberg and Katznelson and
concentrate on proving Theorem~\ref{thm:inf-DHJ}, and will also
follow their next step (Lemma~\ref{lem:sssuffice} below) to reduce
our study to a class of `strongly stationary' probability measures.
Unlike the earlier settings of ergodic Ramsey Theory, this
stationarity condition is not readily described by a collection of
invertible probability-preserving transformations (in particular, it
is instead described in terms of a very large semigroup of highly
non-invertible transformations, which cannot easily be made
invertible by passing to any simple extended system while respecting
the relations of the semigroup). Consequently Furstenberg and
Katznelson must next impose a collection of invertible
transformations `by hand' that describes only a rather weaker
subsemigroup of symmetries, and then bring modifications of their
older ergodic-theoretic techniques~(see, in particular, \cite{Fur81}
and~\cite{FurKat85}) to bear on these.

Here we avoid the introduction of these transformations, and give an
analysis purely in terms of the strong stationarity obtained
initially. This has much more in common with many of the basic
studies of partially exchangeable arrays of random variables,
particularly by Kingman, Hoover, Aldous and Kallenberg: see, for
example, the recent book of Kallenberg~\cite{Kal05} and the
references given there (and also the survey~\cite{Aus--ERH}, which
treats this subject in a very similar formalism to the present paper
and also describes the relations of those developments to other
combinatorial results in extremal hypergraph theory).  Ultimately we
reduce the problem to an application of the `infinitary hypergraph
removal lemma' of Tao~\cite{Tao07} (or, more precisely, of a
cut-down corollary of that lemma first used
in~\cite{Aus--newmultiSzem} to give a very similar new proof of the
Multidimensional Szemer\'edi Theorem).

We will prove our main structural result for strongly stationary
laws as an assertion that any such law is a `factor' of a law with a
particularly simple structure (similar to the structure of the joint
distribution of all the ingredients that are introduced for the
representation of an exchangeable array), and then this structure
will give the reduction to an infinitary removal lemma. This `simple
structure' will be introduced in the definition of `sated' laws in
Section~\ref{sec:insens} below.

We note here that bringing this general program to bear on the task
of proving Theorem~\ref{thm:inf-DHJ} would not have been possible
without a crucial insight that recently emerged from an ongoing open
collaborative project of Bukh, Gowers, Kalai, McCutcheon, O'Donnell,
Solymosi and Tao.  Their ultimate goal was a purely finitary,
combinatorial proof of Theorem~\ref{thm:DHJ}, and as the present
paper neared completion this also seemed to have been realized;
these developments can be followed online~(\cite{Gow(online)}). The
critical observation that we have taken from their work drives our
proof of Theorem~\ref{thm:char} below, but it has been translated
into a very different lexicon from the finitary work
of~\cite{Gow(online)} and we do not attempt to set up a full
dictionary here.  I am also grateful to Tim Gowers and Terence Tao
for helpful suggestions made about earlier drafts of this paper.

\section{Some background from combinatorics}\label{sec:comb-bac}

We write $[N]$ to denote the discrete interval $\{1,2,\ldots,N\}$
and $\cal{P}S$ to denote the power set of $S$.

Most of our work will consider probabilities on product spaces
indexed by the infinite-dimensional combinatorial space $[k]^\o$
introduced above. We will denote the concatenation of two finite
words $u,v \in [k]^\o$ by either $uv$ or $u\oplus v$.  For any fixed
finite $n$ we can define an \textbf{$n$-dimensional subspace} of
$[k]^\o$ to be an injection $\phi:[k]^n \into [k]^\o$ specified as
follows: for some integers $0 = N_0 < N_1 < N_2 < \ldots < N_n$,
nonempty subsets $I_1 \subseteq [N_1]$, $I_2 \subseteq
[N_2]\setminus [N_1]$, \ldots, $I_n \subseteq [N_n]\setminus
[N_{n-1}]$ and fixed words $w_1 \in [k]^{N_1}$, $w_2 \in [k]^{N_2}$,
\ldots, $w_n \in [k]^{N_n}$ we let $\phi(v_1v_2\cdots v_n)$ be the
word in $[k]^\o$ of length $N_n$ such that when $N_i < m \leq
N_{i+1}$ we have
\[\phi(v_1v_2\cdots v_n)_m := \left\{\begin{array}{ll}(w_{i+1})_m&\quad\quad\hbox{if }m\in \{N_i+1,N_i+2,\ldots,N_{i+1}\}\setminus I_{i+1}\\ v_i&\quad\quad\hbox{if }m\in I_{i+1}.\end{array}\right.\]

In these terms a combinatorial line is simply a $1$-dimensional
combinatorial subspace.

Similarly, an \textbf{infinite-dimensional subspace} (or often just
\textbf{subspace}) of $[k]^\o$ is an injection $\phi:[k]^\o \into
[k]^\o$ specified by the above rule for some infinite sequence $0 =
N_0 < N_1 < N_2 < \ldots$ and nonempty $I_{i+1} \subseteq
[N_{i+1}]\setminus [N_i]$.  It is clear that the collection of all
subspaces of $[k]^\o$ forms a semigroup under composition.

Finally, let us define \textbf{letter-replacement maps}: give $i\in
[k]$ and $e \subseteq [k]$, for each $N\geq 1$ we define
$r^N_{e,i}:[k]^N \to [k]^N$ by
\[r^N_{e,i}(w)_m := \left\{\begin{array}{ll}i&\quad\quad\hbox{if }w_m \in e\\ w_m&\quad\quad\hbox{if }w_m \in [k]\setminus e\end{array}\right.\]
for $m\leq N$, and let \[r_{e,i} := \bigcup_{N\geq
1}r^N_{e,i}:[k]^\o\to [k]^\o\] (so clearly $r_{e,i}$ actually takes
values in the subset $([k]\setminus (e\setminus\{i\}))^{[k]^\o}$.

\section{Some background from probability}

Throughout this paper $(X,\S)$ will denote a standard Borel
measurable space. We shall write $(X^I,\S^{\otimes I})$ for the
usual product measurable structure indexed by a set $I$ and
$\mu^{\otimes I}$ for the product of a probability measure $\mu$ on
$(X,\S)$. Given a measurable map $\phi:(X,\S)\to (Y,\Phi)$ to
another standard Borel space, we shall write $\phi_\#\mu$ for the
resulting pushforward probability measure on $(Y,\Phi)$.  We will
generally use $\pi_J$ to denote any coordinate projection from a
product space $X^I$ onto its factor $X^J$ for any $I \supseteq J$,
and will shorten $\pi_{\{j\}}$ to $\pi_j$.

Most of our interest will be in probability measures on the product
spaces $(X^{[k]^\o},\S^{\otimes [k]^\o})$ for various standard Borel
spaces $(X,\S)$.  In this paper we will simply refer to these as
\textbf{laws}, in view of their interpretation as the joint laws of
$X$-valued stochastic processes indexed by $[k]^\o$. Let us note
here that Theorem~\ref{thm:inf-DHJ} is clearly equivalent to the
following superficially more general result, whose formulation will
be more convenient for our proof.

\begin{thm}\label{thm:inf-DHJ2}
For any $\delta > 0$, if $(X,\S)$ is a standard Borel space, $\mu$
is a law on $(X^{[k]^\o},\S^{\otimes [k]^\o})$ and $A \in \S$ is such
that $\mu(\pi_w^{-1}(A))\geq \delta$ for every $w \in [k]^\o$ then
there are some $m\geq 1$ and a combinatorial line $\ell:[k]\into
[k]^m$ such that
\[\mu\Big(\bigcap_{i=1}^k\pi_{\ell(i)}^{-1}(A)\Big) > 0.\]
\end{thm}

If $\mu$ is a law and $\phi:[k]^\o \into [k]^\o$ is a subspace, then
the projected law $(\pi_{\rm{image}(\phi)})_\#\mu$ on
$X^{\rm{image}(\phi)}$ can be canonically identified with
another law on $X^{[k]^\o}$, simply because $\phi$ itself gives an
identification of $[k]^\o$ with $\rm{image}(\phi)$. In this case we
will write $\phi^\ast\mu$ for this new law on $X^{[k]^\o}$.

Borrowing some notation from ergodic theory, a \textbf{factor} of a
law $\mu$ on $(X^{[k]^\o},\S^{[k]^\o})$ will be a Borel map
$\phi:(X,\S)\to (Y,\Phi)$ to some other standard Borel space
$(Y,\Phi)$.  To such a map we can associate its inverse-image
$\s$-subalgebra $\phi^{-1}(\Phi)\leq \S$, and it is standard that in
the category of Borel spaces, given a Borel probability measure on
$\S$ any $\s$-subalgebra of $\S$ agrees with the inverse-image
$\s$-subalgebra of some factor $\phi$ up to modifying by negligible
sets (see, for example, Chapter 2 of Glasner~\cite{Gla03}).  To such
a map $\phi$ we associate the map $\phi^{[k]^\o}:X^{[k]^\o}\to X^{[k]^\o}$
corresponding to the coordinate-wise action of $\phi$, and will
refer to $(\phi^{[k]^\o})_\#\mu$ as the associated \textbf{factor
law} of $\mu$.  In the opposite direction, if $\mu$ arises from a
factor of some `larger' law $\l$ via the factor $\phi$ then we will
refer to $\l$ as an \textbf{extension of $\mu$ through $\phi$}.

An \textbf{inverse system of laws} comprises an inverse system of
standard Borel spaces
\[\ldots \stackrel{\psi^{(m+2)}_{(m+1)}}{\longrightarrow}(X_{(m+1)},\S_{(m+1)})\stackrel{\psi^{(m+1)}_{(m)}}{\longrightarrow}(X_{(m)},\S_{(m)})\stackrel{\psi^{(m)}_{(m-1)}}{\longrightarrow}\ldots\stackrel{\psi^{(1)}_{(0)}}{\longrightarrow}(X_{(0)},\S_{(0)})\]
together with a sequence of laws $\mu_{(m)}$ on
$(X_{(m)}^{[k]^\o},\S_{(m)}^{[k]^\o})$ such that
$((\psi^{(m+1)}_{(m)})^{[k]^\o})_\#\mu_{(m+1)} = \mu_{(m)}$ for
every $m$. In this case we will define
\[\psi^{(m)}_{(k)}:= \psi^{(k+1)}_{(k)}\circ\psi^{(k+2)}_{(k+1)}\circ\cdots\circ\psi^{(m)}_{(m-1)}\]
for $k\leq m$, and will sometimes write instead
\[\ldots \stackrel{\psi^{(m+2)}_{(m+1)}}{\longrightarrow}(X^{[k]^\o}_{(m+1)},\S^{\otimes[k]^\o}_{(m+1)},\mu_{(m+1)})\stackrel{\psi^{(m+1)}_{(m)}}{\longrightarrow}(X^{[k]^\o}_{(m)},\S^{\otimes [k]^\o}_{(m)},\mu_{(m)})\stackrel{\psi^{(m)}_{(m-1)}}{\longrightarrow}\ldots\]
as a shorthand to denote this overall situation.

Given an inverse sequence as above, then exactly as in standard
ergodic theory (see, for example, Examples 6.3 of
Glasner~\cite{Gla03}) we can construct an \textbf{inverse limit} in
the form of a standard Borel space $(X_{(\infty)},\S_{(\infty)})$, a
law $\mu_{(\infty)}$ on
$(X^{[k]^\o}_{(\infty)},\S^{\otimes[k]^\o}_{(\infty)})$ and a family
of factors $\psi_{(m)}:X_{(\infty)}\to X_{(m)}$ such that
$\psi_{(k)} = \psi^{(m)}_{(k)}\circ\psi_{(m)}$ for all $k < m$ and
$((\psi_{(k)})^{[k]^\o})_\#\mu_{(\infty)} = \mu_{(m)}$ for every
$m$.  We will use this construction later in the paper.

Related to the notion of a factor is that of a `coupling': given
laws $\mu$ and $\nu$ on $(X^{[k]^\o},\S^{\otimes[k]^\o})$ and
$(Y^{[k]^\o},\Phi^{\otimes[k]^\o})$ respectively, a
\textbf{coupling} of $\mu$ and $\nu$ is a law $\l$ on $((X\times
Y)^{[k]^\o},(\S\otimes \Phi)^{\otimes[k]^\o})$ whose coordinate
projections onto $(X^{[k]^\o},\S^{\otimes[k]^\o})$ and
$(Y^{[k]^\o},\Phi^{\otimes[k]^\o})$ are $\mu$ and $\nu$
respectively.  This definition generalizes to couplings of larger
collections of laws in the obvious way.  We will also have need for
a topology on couplings, set up exactly analogously with the
`joining topology' of ergodic theory: quite generally, given a
countable collection of standard Borel probability spaces
$(X_i,\S_i,\mu_i)_{i\in I}$, the space $C$ of all couplings of the
$\mu_i$ on the product standard Borel space $\big(\prod_{i\in
I}X_i,\bigotimes_{i\in I}\S_i\big)$ is endowed with the weakest
topology with respect to which all the evaluation maps
\[\l \mapsto \int_{\prod_{i\in I}X_i}\prod_{i \in F}f_i\circ\pi_i\,\d\l\]
for collections $f_i \in L^\infty(\mu_i)$ indexed by finite subsets
$F \subseteq I$ are continuous.  Just as for joinings of
probability-preserving systems (as discussed in Chapter 6 of
Glasner~\cite{Gla03}), the restriction here to couplings of fixed
one-dimensional marginals (rather than arbitrary probability
measures on the product space) gives that this is a compact topology
on $C$.

\section{Strongly stationary laws}

We now introduce the special class of laws that will concern us
through most of this paper.  These are distinguished by satisfying a
kind of `self-similarity' in terms of the structure of the index set
$[k]^\o$.

\begin{dfn}[Strong stationarity]
A law $\mu$ on $(X^{[k]^\o},\S^{[k]^\o})$ is \textbf{strongly
stationary} (\textbf{s.s.}) if $\phi^\ast\mu = \mu$ for every
subspace $\phi:[k]^\o\into [k]^\o$.
\end{dfn}

This can be thought of as the analog appropriate to the present
setting of the exchangeability of a family of random variables (or,
equivalently, their joint distribution) under an index-set-permuting
action of some countable group: see, for example, Section 2.2
of~\cite{Aus--ERH}, where this abstract definition is set up before
being applied to exchangeable arrays (or `exchangeable random
hypergraphs', as they are formulated there).

Indeed, the only real difference between the settings of that paper
and this is that here our notion of strong stationarity refers to a
semigroup of noninvertible self-maps of the underlying index set,
for which it seems difficult to find any `invertible model'.
Furstenberg and Katznelson meet the same difficulty in their
original work, and circumvent it by relying instead only on a weaker
symmetry to which they can associate (using a highly arbitrary
selection procedure) a collection of invertible
probability-preserving transformations. By contrast, we will find
that this noninvertibility is of no consequence for our approach
below.

Let us next recall Furstenberg and Katznelson's reduction to the
case of s.s. laws, contained in Sections 2.3 and 2.4
of~\cite{FurKat91}.

\begin{lem}\label{lem:sssuffice}
If Theorem~\ref{thm:inf-DHJ2} holds for all s.s. laws for every
$\delta > 0$ then it holds for all laws for every $\delta > 0$.
\end{lem}

\textbf{Proof}\quad We only sketch the argument, referring the
reader to~\cite{FurKat91} for the details.  First note that given a law $\mu$ on $(X^{[k]^\o},\S^{\otimes [k]^\o})$ for a general standard Borel space $(X,\S)$ and a subset $A \in \S$, the joint distribution under $\mu$ of the sets $\pi_w^{-1}(A)$, $w \in [k]^\o$, defines a law on $\{0,1\}^{[k]^\o}$, and the desired conclusion depends only on this factor law, so we may assume that $X = \{0,1\}$.

Now the point is that by
applying the Carlson-Simpson Theorem~\cite{CarSim84} (see also~\cite{FurKat89}) to arbitrarily
fine finite coverings of the finite-dimensional spaces of
probability distributions on $\{0,1\}^{[k]^n}$ for increasingly
large $n$, we obtain a subspace $\psi:[k]^\o \into [k]^\o$ and an
infinite word $w \in [k]^\bbN$ such that the restricted laws
\[\psi(w|_{[m]}\oplus \cdot)^\ast\mu\]
converge to a strongly stationary law as $m\to\infty$, and since
\emph{all} one-dimensional marginals of the input law gave
probability at least $\delta$ to $\{1\}$, the same is true of the
limit.  Finally, the condition
\[\mu\{\bf{x} \in \{0,1\}^{[k]^\o}:\ x_{\ell(i)} = 1\ \forall i\leq k\} > 0\]
is also finite-dimensional and open for this topology on the space
of finite-dimensional distributions, so if it holds for the limit
measure it must also hold somewhere for the original measure. \qed

\begin{dfn}
If the law $\mu$ is s.s. then in particular all the one-dimensional
marginals $(\pi_w)_\#\mu$ for $w \in [k]^\o$ are the same and all
the $k$-dimensional marginals
\[(\pi_{\ell(1)},\pi_{\ell(2)},\ldots,\pi_{\ell(k)})_\#\mu\]
for $\ell$ a line in $[k]^\o$ the same.  We will refer to these as
the \textbf{point-marginal} and \textbf{line-marginal} of $\mu$ and
will often denote them by $\mu^\circ$ and $\mu^\line$ respectively.
\end{dfn}

\section{Partially insensitive and sated
laws}\label{sec:insens}

\begin{dfn}[Partially insensitive $\s$-algebras]
For any nonempty $e \subseteq [k]$ and a s.s. law $\mu$ on
$(X^{[k]^\o},\S^{\otimes [k]^\o})$ the \textbf{$e$-insensitive
$\s$-algebra} is the $\s$-subalgebra $\Phi_e \leq \S$ defined by
\[\Phi_e := \{A \in\S:\ 1_A(x_{\ell(i)}) = 1_A(x_{\ell(j)})\ \forall i,j\in e\ \hbox{for}\ \mu\hbox{-a.e.}\ (x_w)_w\in X^{[k]^\o}\}\]
(note that any choice of line $\ell$ will do here, owing to the assumption that $\mu$ is s.s.). The $e$-insensitive $\s$-algebras for different sets $e$ are
together referred to as the \textbf{partially insensitive
$\s$-algebra}.  A measurable function $f$ on $X$ is
\textbf{$e$-insensitive} if it is $\Phi_e$-measurable.

The law $\mu$ is itself \textbf{$e$-insensitive} if $\Phi_e = \S$,
that is if $x_{\ell(i)} = x_{\ell(j)}$ for every $i,j \in e$ for
$\mu$-a.e. $(x_w)_w\in X^{[k]^\o}$.
\end{dfn}

We now also construct a larger collection of $\s$-algebras from the
above, but first must set up some additional notation. These next
$\s$-algebras will be indexed by \textbf{up-sets} in
$\binom{[k]}{\geq 2}$: that is, families $\I\subseteq
\binom{[k]}{\geq 2}$ such that if $u \in \I$ and $[k]\supseteq
v\supseteq u$ then also $v \in \I$. For example, given $e \subseteq
[k]$ we write $\langle e\rangle := \{u \in\binom{[k]}{\geq 2}:\
u\supseteq e\}$ (note the non-standard feature of our notation that
$e \in \langle e\rangle$ if and only if $|e| \geq 2$): up-sets of
this form are \textbf{principal}.  We will abbreviate
$\langle\{i\}\rangle$ to $\langle i\rangle$.

In general, for any up-set $\I\subseteq \binom{[k]}{\geq 2}$ we let
$\Phi_\I:= \bigvee_{e \in \I}\Phi_e$.  It is clear from the above
definition that if $e \subseteq e'$ then $\Phi_e \supseteq
\Phi_{e'}$, so we have $\Phi_e = \Phi_{\langle e\rangle}$.

It is also immediate from the above definition that for any s.s. law
$\mu$, $e \in \binom{[k]}{\geq 2}$ and $i,j\in e$ the
$\s$-subalgebras $\pi_i^{-1}(\Phi_e)$ and $\pi_j^{-1}(\Phi_e)$ of
$\S^{\otimes k}$ are equal up to $\mu^\line$-negligible sets, and so
we can make the following definition.

\begin{dfn}[Oblique copies]
For each $e \subseteq [k]$ we refer to the common
$\mu^\line$-completion of the $\s$-subalgebra $\pi_i^{-1}(\Phi_e)$,
$i \in e$, as the \textbf{oblique copy} of $\Phi_e$, and denote it
by $\Phi^\dag_e$.  More generally we shall refer to $\s$-algebras
formed from the oblique copies by repeatedly applying $\cap$ and $\vee$ as \textbf{oblique $\s$-algebras}, and if $\I \subseteq \binom{[k]}{\geq 2}$ is any up-set then we let $\Phi^\dag_\I := \bigvee_{e \in \I}\Phi^\dag_e$.
\end{dfn}

Clearly if a law is $e$-insensitive for some $e$ this amounts to a
nontrivial simplification of its structure.  In general we will
analyze an arbitrary law in terms of its possible couplings to
insensitive laws through the following definition.

\begin{dfn}[Sated laws]
For a nonempty up-set $\I \subseteq \binom{[k]}{\geq 2}$ and a s.s.
law $\mu$ on $(X^{[k]^\o},\S^{\otimes [k]^\o})$ with partially
insensitive $\s$-algebras $\Phi_e$, $\mu$ is \textbf{$\I$-sated} if
for any s.s. extension $\tilde{\mu}$ of $\mu$ the factor
$\pi:\tilde{x} \mapsto x$ and the $\s$-subalgebra $\tilde{\Phi}_\I$
are relatively independent under $\tilde{\mu}^\circ$ over the
$\s$-subalgebra $\pi^{-1}(\Phi_\I)$.

The law $\mu$ is \textbf{fully sated} if it is $\I$-sated for every
such $\I$.
\end{dfn}

Clearly not all laws are sated, but it turns out that we can recover
the advantage of working with a sated law by passing to an
extension.  The following theorem is closely analogous to a similar
`satedness' result to appear in~\cite{Aus--lindeppleasant}, and is
also closely related to older results on `pleasant' and
`isotropized' extensions of probability-preserving systems
in~\cite{Aus--newmultiSzem,Aus--nonconv}.

\begin{thm}[Sated extension]\label{thm:satedsexist}
Every s.s. law has a fully sated s.s. extension.
\end{thm}

In light of this it will suffice to prove Theorem~\ref{thm:inf-DHJ2}
for fully sated s.s. laws $\mu$.  We will finish this section by
proving Theorem~\ref{thm:satedsexist}, and then in the next section
we will derive some useful consequences of full satedness for the
structure of $\mu^\line$ before using these to complete the
reduction of Theorem~\ref{thm:inf-DHJ2} to an infinitary removal
lemma in Section~\ref{sec:final}.

\begin{lem}[Partially sated extension]
For any up-set $\I\subseteq \binom{[k]}{\geq 2}$, every s.s. law
$\mu$ has an s.s. extension that is $\I$-sated.
\end{lem}

\textbf{Proof}\quad This proceeds by an infinitary `energy
increment' argument: we build a tower of extensions of $\mu$ each
`closer' to $\I$-satedness than its predecessor and so that the
resulting inverse limit is exactly $\I$-sated.

Let $(f_r)_{r\geq 1}$ be a countable subset of the $L^\infty$-unit
ball $\{f\in L^\infty(\mu^\circ):\ \|f\|_\infty\leq 1\}$ that is
dense in this ball for the $L^2$-norm, and let $(r_i)_{i\geq 1}$ be
a member of $\bbN^\bbN$ in which every non-negative integer appears
infinitely often.

We will now construct an inverse sequence
\[\ldots\stackrel{\psi^{(m+2)}_{(m+1)}}{\longrightarrow}(X_{(m+1)}^{[k]^\o},\S^{\otimes [k]^\o}_{(m+1)},\mu_{(m+1)})\stackrel{\psi^{(m+1)}_{(m)}}{\longrightarrow} (X^{[k]^\o}_{(m)},\S^{\otimes [k]^\o}_{(m)},\mu_{(m)})\stackrel{\psi^{(m)}_{(m-1)}}{\longrightarrow}\ldots\]
starting from $(X_{(0)},\S_{(0)}) = (X,\S)$ and $\mu_{(0)} = \mu$
such that each $(X_{(m+1)}^{[k]^\o},\S^{\otimes
[k]^\o}_{(m+1)},\mu_{(m+1)})$ is obtained by coupling to
$(X_{(m)}^{[k]^\o},\S^{\otimes [k]^\o}_{(m)},\mu_{(m)})$ a new law
$\mu'$ on some $((X')^{[k]^\o},(\S')^{\otimes [k]^\o})$ such that
for this new law we have $\S' = \Phi'_\I$.

Suppose that we have already obtained $(X_{(m)}^{[k]^\o},\S^{\otimes
[k]^\o}_{(m)},\mu_{(m)})$ for $0 \leq m \leq m_1$.  We consider two
separate cases:
\begin{itemize}
\item If there is some further extension
\[\pi:(\tilde{X}^{[k]^\o},\tilde{\S}^{\otimes[k]^\o},\tilde{\mu})\to (X^{[k]^\o}_{(m_1)},\S_{(m_1)}^{\otimes [k]^\o},\mu_{(m_1)})\]
such that
\[\|\sfE_{\tilde{\mu}^\circ}(f_{r_{m_1}}\circ \psi^{(m_1)}_{(0)}\circ\pi\,|\,\tilde{\Phi}_\I)\|_2^2 > \|\sfE_{\mu_{(m_1)}^\circ}(f_{r_{m_1}}\circ \psi^{(m_1)}_{(0)}\,|\,\Phi_{(m_1),\I})\|_2^2 + 2^{-m_1},\]
then choose a particular such extension such that the increase
\[\|\sfE_{\tilde{\mu}^\circ}(f_{r_{m_1}}\circ \psi^{(m_1)}_{(0)}\circ\pi\,|\,\tilde{\Phi}_\I)\|_2^2 - \|\sfE_{\mu_{(m_1)}^\circ}(f_{r_{m_1}}\circ \psi^{(m_1)}_{(0)}\,|\,\Phi_{(m_1),\I})\|_2^2\]
is at least half its supremal possible value over such extensions.
Now by restricting to the possibly smaller extension of
$(X^{[k]^\o}_{(m_1)},\S_{(m_1)}^{\otimes [k]^\o},\mu_{(m_1)})$ given
by replacing $(\tilde{X},\tilde{\S})$ with its factor generated by
$\pi$ and the $\s$-algebra $\tilde{\Phi}_\I$, we may assume that
$\tilde{\mu}$ is itself obtained as a coupling of $\mu_{(m_1)}$ to a
law $\mu'$ for which the $\s$-algebra $\Phi'_\I$ is full, and now we
let $(X_{(m_1+1)},\S_{(m_1+1)}) := (\tilde{X},\tilde{\S})$,
$\mu_{(m_1 + 1)} := \tilde{\mu}$ and $\psi^{(m_1 + 1)}_{(m_1)} :=
\pi$.
\item If, on the other hand, for every further extension $\pi$ as above we
have
\[\|\sfE_{\tilde{\mu}^\circ}(f_{r_{m_1}}\circ
\psi^{(m_1)}_{(0)}\circ\pi\,|\,\tilde{\Phi}_\I)\|_2^2 \leq
\|\sfE_{\mu_{(m_1)}^\circ}(f_{r_{m_1}}\circ
\psi^{(m_1)}_{(0)}\,|\,\Phi_{(m_1),\I})\|_2^2 + 2^{-m_1},\] then we
simply set $\psi^{(m_1+1)}_{(m_1)}:= \id_{X_{(m_1)}}$.
\end{itemize}

Finally, let $(X_{(\infty)},\S_{(\infty)},\mu_{(\infty)}^\circ)$ be
the inverse limit probability space of
\[\ldots\stackrel{\psi^{(m+2)}_{(m+1)}}{\longrightarrow}(X_{(m+1)},\S_{(m+1)},\mu^\circ_{(m+1)})\stackrel{\psi^{(m+1)}_{(m)}}{\longrightarrow} (X_{(m)},\S_{(m)},\mu_{(m)}^\circ)\stackrel{\psi^{(m)}_{(m-1)}}{\longrightarrow}\ldots,\]
$\mu_{(\infty)}$ the inverse limit of the measures $\mu_{(m)}$ and
$\psi_{(m)}:X_{(\infty)}\to X_{(m)}$ the resulting factor maps. It
is clear from the above construction that the whole $\s$-algebra
$\S_{(\infty)}$ is generated up to $\mu_{(\infty)}^\circ$-negligible
sets by $\Phi_{(\infty),\I}$ and $\psi_{(0)}$, since
$\Phi_{(\infty),\I}$ contains every
$\psi_{(m)}^{-1}(\Phi_{(m),\I})$. To show that $\mu_{(\infty)}$ is
$\I$-sated, let $\tilde{\mu}$ under $\pi:\tilde{X}\to X_{(\infty)}$
be any further extension of $\mu_{(\infty)}$, and suppose that $f
\in L^\infty(\mu_{(\infty)})$. We will complete the proof by showing
that
\[\sfE_{\tilde{\mu}^\circ}(f\circ\pi\,|\,\tilde{\Phi}_\I) =
\sfE_{\mu_{(\infty)}^\circ}(f\,|\,\Phi_{(\infty),\I})\circ\pi.\]

By construction, this $f$ may be approximated arbitrarily well in
$L^2(\mu^\circ_{(\infty)})$ by finite sums of the form $\sum_p
g_p\cdot h_p$ with $g_p$ being bounded and
$\Phi_{(\infty),\I}$-measurable and $h_p$ being bounded and
$\psi_{(0)}$-measurable, and now by density we may also restrict to
using $h_p$ that are each a scalar multiple of some
$f_{r_p}\circ\psi_{(0)}$, so by continuity and multilinearity it
suffices to prove the above equality for one such product $g\cdot
(f_r\circ\psi_{(0)})$.  Since $g$ is $\Phi_{(\infty),\I}$-measurable
and $\tilde{\Phi}_\I \supseteq \pi^{-1}(\Phi_{(\infty),\I})$, it will now be sufficient to show that
\[\sfE_{\tilde{\mu}^\circ}(f_r\circ\psi_{(0)}\circ\pi\,|\,\tilde{\Phi}_\I) =
\sfE_{\mu_{(\infty)}^\circ}(f_r\circ\psi_{(0)}\,|\,\Phi_{(\infty),\I})\circ\pi,\]
and this in turn will follow if we only show that
\[\|\sfE_{\tilde{\mu}^\circ}(f_r\circ\psi_{(0)}\circ\pi\,|\,\tilde{\Phi}_\I)\|_2^2 =
\|\sfE_{\mu_{(\infty)}^\circ}(f_r\circ\psi_{(0)}\,|\,\Phi_{(\infty),\I})\|_2^2.\]

Now, by the martingale convergence theorem we have
\[\|\sfE_{\mu^\circ_{(m)}}(f_r\circ\psi^{(m)}_{(0)}\,|\,\Phi_{(m),\I})\|_2^2 \uparrow \|\sfE_{\mu^\circ_{(\infty)}}(f_r\circ\psi_{(0)}\,|\,\Phi_{(\infty),\I})\|_2^2\]
as $m\to\infty$.  It follows that if
\[\|\sfE_{\tilde{\mu}^\circ}(f_r\circ\psi_{(0)}\circ\pi\,|\,\tilde{\Phi}_\I)\|_2^2
>
\|\sfE_{\mu^\circ_{(\infty)}}(f_r\circ\psi_{(0)}\,|\,\Phi_{(\infty),\I})\|_2^2\]
then for some sufficiently large $m$ we would have $r_m = r$ (since
each integer appears infinitely often as some $r_m$) but
\begin{eqnarray*}
&&\|\sfE_{\mu^\circ_{(m+1)}}(f_{r_m}\circ\psi^{(m+1)}_{(0)}\,|\,\Phi_{(m+1),\I})\|_2^2
-
\|\sfE_{\mu^\circ_{(m)}}(f_r\circ\psi^{(m)}_{(0)}\,|\,\Phi_{(m),\I})\|_2^2\\
&& \leq
\|\sfE_{\mu^\circ_{(\infty)}}(f_r\circ\psi_{(0)}\,|\,\Phi_{(\infty),\I})\|_2^2
-
\|\sfE_{\mu^\circ_{(m)}}(f_r\circ\psi^{(m)}_{(0)}\,|\,\Phi_{(m),\I})\|_2^2\\
&& <
\frac{1}{2}\Big(\|\sfE_{\tilde{\mu}^\circ}(f_r\circ\psi^{(\infty)}_{(0)}\circ\pi\,|\,\tilde{\Phi}_\I)\|_2^2
-
\|\sfE_{\mu^\circ_{(m)}}(f\circ\psi^{(m)}_{(0)}\,|\,\Phi_{(m),\I})\|_2^2\Big)
\end{eqnarray*}
and also
\[\|\sfE_{\tilde{\mu}^\circ}(f_r\circ\psi^{(\infty)}_{(0)}\circ\pi\,|\,\tilde{\Phi}_\I)\|_2^2\geq
\|\sfE_{\mu^\circ_{(m)}}(f\circ\psi^{(m)}_{(0)}\,|\,\Phi_{(m),\I})\|_2^2
+ 2^{-m}\] so contradicting our choice of $\mu_{(m+1)}$ in the first
alternative in our construction above. This contradiction shows that
we must actually have the equality of $L^2$-norms asserted above, as
required. \qed

\textbf{Proof of Theorem~\ref{thm:satedsexist}}\quad Pick a sequence
of up-sets $(\I_m)_{m\geq 1}$ in which each possible up-set appears
infinitely often. Now by repeatedly implementing the preceding lemma
we can form another tower of extensions
\[\ldots\to(X_{(m+1)}^{[k]^\o},\S^{\otimes [k]^\o}_{(m+1)},\mu_{(m+1)})\to (X^{[k]^\o}_{(m)},\S^{\otimes [k]^\o}_{(m)},\mu_{(m)})\to\ldots\]
above $(X^{[k]^\o},\S^{\otimes [k]^\o},\mu)$ in which every
$\mu^{(m)}$ is $\I_m$-sated.  It is now an immediate check that the
resulting inverse limit $(\tilde{X}^{[k]^\o},\tilde{\S}^{\otimes
[k]^\o},\tilde{\mu})$ is fully sated. \qed

\section{The structure of sated laws}

Having proved the existence of sated extensions, we will now show
how the structure of $\mu$ (and particularly of the partially
insensitive $\s$-algebras $\Phi_e$) simplifies for sated systems,
before using these results to prove Theorem~\ref{thm:inf-DHJ2} in
the next section.

First we need the following lemma.

\begin{lem}\label{lem:isotropized}
If $\mu$ is fully sated then for every $i\in e \in \binom{[k]}{\geq
2}$, if $f \in L^\infty(\mu^\circ)$ is $e$-insensitive then
\[\sfE_{\mu^\circ}\Big(f\,\Big|\,\bigvee_{j\in[k]\setminus e}\Phi_{\{i,j\}}\Big) = \sfE_{\mu^\circ}\Big(f\,\Big|\,\bigvee_{j\in[k]\setminus e}\Phi_{e\cup\{j\}}\Big).\]
\end{lem}

\textbf{Proof}\quad Clearly
\[\sfE_{\mu^\circ}\Big(f\,\Big|\,\bigvee_{j\in[k]\setminus e}\Phi_{e\cup\{j\}}\Big)\]
is always $\big(\bigvee_{j\in[k]\setminus
e}\Phi_{\{i,j\}}\big)$-measurable. It will therefore suffice to show
that if $f \in L^\infty(\mu^\circ)$ is $e$-insensitive and
orthogonal to the $\s$-algebra $\bigvee_{j\in [k]\setminus
e}\Phi_{e\cup\{j\}}$ then it is actually orthogonal to
$\bigvee_{j\in[k]\setminus e}\Phi_{\{i,j\}}$.  We prove this by
contradiction, so suppose for one such $f$ that we could find some
bounded functions $h_j$ for $j\in[k]\setminus e$ such that each
$h_j$ is $\Phi_{\{i,j\}}$-measurable and
\[\int_X f\cdot \prod_{j\in [k]\setminus e}h_j\,\d\mu^\circ = \kappa \neq 0.\]
Re-writing this inner product condition at the level of the whole
law $\mu$ it simply reads that
\[\int_{X^{[k]^\o}}f(x_w)\cdot \prod_{j\in [k]\setminus e}h_j(x_w)\,\mu(\d\bf{x})  = \kappa\]
for any fixed $w \in [k]^\o$.  However, now we apply first the
$e$-insensitivity of $f$ to deduce that also
\begin{multline*}
\int_{X^{[k]^\o}}f(x_w)\cdot \prod_{j\in [k]\setminus
e}h_j(x_{r_{e,i}(w)})\,\mu(\d\bf{x})\\
=\int_{X^{[k]^\o}}f(x_{r_{e,i}(w)})\cdot \prod_{j\in [k]\setminus
e}h_j(x_{r_{e,i}(w)})\,\mu(\d\bf{x}) = \kappa
\end{multline*}
(where $r_{e,i}$ is the letter-replacement map defined at the end of Section~\ref{sec:comb-bac}) for every word $w$, and now the $\{i,j\}$-insensitivity of $h_j$ to
deduce that
\begin{multline*}
\int_{X^{[k]^\o}}f(x_w)\cdot \prod_{j\in [k]\setminus
e}h_j(x_{r_{e,j}(w)})\,\mu(\d\bf{x})\\
=\int_{X^{[k]^\o}}f(x_w)\cdot \prod_{j\in [k]\setminus
e}h_j(x_{r_{e,i}(w)})\,\mu(\d\bf{x}) = \kappa.
\end{multline*}
for every word $w$.

It follows that if we define the probability measure $\l$ on
$(X\times X^{[k]\setminus e})^{[k]^\o}$ to be the joint law under
$\mu$ of
\[(x_w)_w \mapsto \big(x_w,(x_{r_{e,j}(w)})_{j\in [k]\setminus e}\big)_w\]
then all of its coordinate projections onto individual copies of $X$
are still just $\mu^\circ$, the projection
\[\pi:\big(y_w,(z_{j,w})_{j\in [k]\setminus e}\big)_w\to (y_w)_w\]
has $\pi_\#\l = \mu$ and the projections
\[\pi_j:\big(y_w,(z_{j,w})_{j\in [k]\setminus e}\big)_w\to (z_{j,w})_w\]
are $\l$-almost surely $(e\cup\{j\})$-insensitive.  Therefore
through the first coordinate projection $\pi$ the law $\l$ defines
an extension of $\mu$, and the above inequality gives a non-zero
inner product under $\l$ for $f$ with some product over $j\in
[k]\setminus e$ of
$(e\cup\{j\})$-insensitive functions, which we can express as
\[\int_{X^{[k]^\o}} (f\circ\pi)\cdot\prod_{j \in [k]\setminus e}(h_j\circ\pi_j)\,\d\l  = \kappa.\]

Now $\l$ may not be stationary, but at least its marginals onto all
individual copies of $X$ in $(X\times X^{[k]\setminus e})^{[k]^\o}$
are equal to $\mu^\circ$.  It follows that we can re-run the appeal
to the Carlson-Simpson Theorem in Lemma~\ref{lem:sssuffice} to
obtain a subspace $\psi:[k]^\o\into [k]^\o$ and an infinite word $w
\in [k]^\bbN$ such that the pulled-back measures
\[\psi(w|_{[m]}\oplus \cdot)^\ast\l\]
converge in the coupling topology on $(X\times X^{[k]\setminus
e})^{[k]^\o}$ (recall that for couplings of fixed marginals this is
compact) to a s.s. measure $\tilde{\mu}$. Since $\mu$ was already
strongly stationary, we must still have $\pi_\#\tilde{\mu}= \mu$,
and by the definition of the coupling topology as the weakest for
which integration of fixed product functions is continuous it
follows that we must still have, firstly, that
\[\int_{X^{[k]^\o}}(f\circ\pi)\cdot\prod_{j \in [k]\setminus e}(h_j\circ\pi_j)\,\d\tilde{\mu}  = \kappa,\]
and secondly that the coordinate projections $\pi_j$ are
$(e\cup\{j\})$-insensitive under $\tilde{\mu}$, since this is
equivalent to the assertion that for any $A \in \S$, $i \in e$ and
line $\ell:[k]\into [k]^\o$ we have
\[\int_{(X\times X^{[k]\setminus e})^{[k]^\o}} 1_A(z_{j,\ell(i)})\cdot 1_{X\setminus A}(z_{j,\ell(j)})\,\tilde{\mu}(\d\bf{z}) = 0\]
and this is clearly closed in the coupling topology.

Therefore we have found a s.s. extension $\tilde{\mu}$ of $\mu$
through some factor map $\xi$ under which
\[\sfE_{\tilde{\mu}^\circ}\Big(f\circ\xi\,\Big|\,\bigvee_{j\in [k]\setminus e}\tilde{\Phi}_{e\cup\{j\}}\Big)\neq 0.\]
By satedness, it follows that in fact
\[\sfE_{\mu^\circ}\Big(f\,\Big|\,\bigvee_{j\in [k]\setminus e}\Phi_{e\cup\{j\}}\Big)\neq 0,\]
contradicting the condition that $f$ be orthogonal to this
$\s$-algebra. \qed

\textbf{Example}\quad The idea behind the above proof may be made
clear by an explication of the special case $k=3$, $i=2$ and $e =
\{1,2\}$. In this case we wish to prove that if $\mu$ is a fully
sated s.s. law on $X^{[3]^\o}$ and $f \in L^\infty(\mu^\circ)$ is
$\{1,2\}$-insensitive then
\[\sfE_{\mu^\circ}(f\,|\,\Phi_{\{2,3\}}) = \sfE_{\mu^\circ}(f\,|\,\Phi_{\{1,2,3\}}),\]
and so we suppose that the right-hand side above is zero and prove
that the left-hand side is also zero.  Arguing by contradiction, we
suppose otherwise and in this case let
\[h := \sfE_{\mu^\circ}(f\,|\,\Phi_{\{2,3\}}).\]
As a $\Phi_{\{2,3\}}$-measurable function, this $h$ must be
$\{2,3\}$-insensitive, and so the condition $h \neq 0$ implies (from
the definition of $h$) that
\[\int_X h^2\,\d\mu^\circ = \int_X fh\,\d\mu^\circ \neq 0,\]
so we have obtained a nontrivial inner product between the
$\{1,2\}$-insensitive function $f$ and the $\{2,3\}$-insensitive
function $h$.  We wish to deduce from this that $f$ actually has a
non-zero inner product with some $\{1,2,3\}$-insensitive function.
For a general s.s. law $\mu$ this does not follow, but using in turn
the strong stationarity of $\mu$, the $\{1,2\}$-insensitivity of $f$
and then the $\{2,3\}$-insensitivity of $h$ we can write
\begin{eqnarray*}
0\neq \int_X fh\,\d\mu^\circ = \int_{X^{[k]^\o}}f(x_w)h(x_w)\,\mu(\d\bf{x}) &=&
\int_{X^{[k]^\o}}f(x_{r_{1,2}(w)})h(x_{r_{1,2}(w)})\,\mu(\d\bf{x})\\
&=& \int_{X^{[k]^\o}}f(x_w)h(x_{r_{1,2}(w)})\,\mu(\d\bf{x})\\
&=&
\int_{X^{[k]^\o}}f(x_w)h(x_{r_{3,2}(r_{1,2}(w))})\,\mu(\d\bf{x})\\
&=& \int_{X^{[k]^\o}}f(x_w)h(x_{222\cdots 2})\,\mu(\d\bf{x})
\end{eqnarray*}
for any $w \in [k]^\o$, where $22\cdots 2$ has the same word-length
as $w$, and this implicitly defines a non-trivial coupling of $\mu$
to a process that is indexed by $\{2\}^\o \subset [3]^\o$ and which
can now be re-interpreted simply as a $\{1,2,3\}$-insensitive law.
Applying the Carlson-Simpson Theorem to construct from this a
similarly nontrivial coupling that is itself s.s. gives a
contradiction with the additional condition that $\mu$ be $\I$-sated
for $\I = \{1,2,3\}$. \fin

The usefulness of satedness for proving Theorem~\ref{thm:inf-DHJ2}
will rest on the following property.

\begin{thm}\label{thm:char}
If $e \subseteq [k]$ is nonempty, $\mu$ is fully sated and $f_i \in
L^\infty(\mu^\circ)$ for $i\in e$ then
\[\int_{X^k}\prod_{i\in e} f_i\circ\pi_i\,\d\mu^\line = \int_{X^k}\prod_{i\in e} \sfE_{\mu^\circ}\Big(f_i\,\Big|\,\bigvee_{l \in e\setminus \{i\}}\Phi_{\{i,l\}}\Big)\circ\pi_i\,\d\mu^\line.\]
\end{thm}

\textbf{Proof}\quad We will prove this by contradiction, assuming
that the desired equality fails for some choice of $f_i \in
L^\infty(\mu^\circ)$ and constructing from this an extension of
$\mu$ witnessing that it is not sated.  For convenience let us
temporarily write $\Xi_i := \bigvee_{l\in
e\setminus\{i\}}\Phi_{\{i,l\}}$ (so $\Xi_i = \Phi_{\langle
i\rangle}$ when $e = [k]$).

Indeed, given such $f_i$ we can write
\begin{multline*}
\int_{X^k}\prod_{i\in e} f_i\circ\pi_i\,\d\mu^\line -
\int_{X^k}\prod_{i\in e} \sfE_{\mu^\circ}(f_i\,|\,\Xi_i)\circ\pi_i\,\d\mu^\line\\
= \sum_{j\in e}\int_{X^k} \Big(\prod_{i \in e,\,i <
j}f_i\circ\pi_i\Big)\cdot(f_j\circ\pi_j -
\sfE_{\mu^\circ}(f_j\,|\,\Xi_j)\circ\pi_j)\cdot\Big(\prod_{i \in
e,\,i >
j}\sfE_{\mu^\circ}(f_i\,|\,\Xi_i)\circ\pi_i\Big)\,\d\mu^\line,
\end{multline*}
and so if this is nonzero then there is some choice of $j\in e$ for
which
\[\int_{X^k} \Big(\prod_{i \in e,\,i <
j}f_i\circ\pi_i\Big)\cdot(f_j\circ\pi_j -
\sfE_{\mu^\circ}(f_j\,|\,\Xi_j)\circ\pi_j)\cdot\Big(\prod_{i \in
e,\,i >
j}\sfE_{\mu^\circ}(f_i\,|\,\Xi_i)\circ\pi_i\Big)\,\d\mu^\line \neq
0.\]

Now for each $i\in e\setminus \{j\}$ recall that
$r_{j,i}:[k]^{\o}\to[k]^{\o}$ is the letter-replacement map defined
by
\[(r_{j,i}(w))_m := \left\{\begin{array}{ll}i&\quad\quad\hbox{if }w_m = j\\ w_m&\quad\quad\hbox{else.}\end{array}\right.\]
In view of the strong stationarity of $\mu$, the above inequality
implies that
\[\int_{X^{[k]^\o}} \Big(\prod_{i \in e,\,i < j}f_i\circ\pi_{r_{j,i}(w)}\Big)\cdot (f_j\circ\pi_w -
\sfE_{\mu^\circ}(f_j\,|\,\Xi_j)\circ\pi_w)\cdot \Big(\prod_{i \in
e,\,i
> j}\sfE_{\mu^\circ}(f_i\,|\,\Xi_i)\circ\pi_{r_{j,i}(w)}\Big)\,\d\mu \neq 0\]
for any $w \in [k]^\o$ such that $w^{-1}\{j\} \neq \emptyset$, since
then the points $r_{j,s}(w)$ for $s=1,2,\ldots,k$ form a
combinatorial line.

It follows that if we define the probability measure $\l$ on
$(X\times X^{e\setminus \{j\}})^{[k]^\o}$ to be the joint law under
$\mu$ of
\[(x_w)_w \mapsto \big(x_w,(x_{r_{j,i}(w)})_{i \in e,\,i < j},(x_{r_{j,i}(w)})_{i \in e,\,i > j}\big)_w\]
then all of its coordinate projections onto individual copies of $X$
are still just $\mu^\circ$, the projection
\[\pi:\big(y_w,(z_{i,w})_{i < j},(z_{i,w})_{i > j}\big)_w\to (y_w)_w\]
has $\pi_\#\l = \mu$ and the projection
\[\pi_{i_0}:\big(y_w,(z_{i,w})_{i \in e,\,i < j},(z_{i,w})_{i \in e,\,i > j}\big)_w\to (z_{i_0,w})_w\]
for $i_0 \in e\setminus\{j\}$ is $\l$-almost surely
$\{i_0,j\}$-insensitive. Therefore through the first coordinate
projection $\pi$ the law $\l$ defines an extension of $\mu$ (not
necessarily s.s.), and the above inequality gives a fixed non-zero
inner product under $\l$ for the function
\[f_j\circ\pi_w -
\sfE_{\mu^\circ}(f_j\,|\,\Xi_j)\circ\pi_w\]
with some product over $i\in e\setminus \{j\}$ of $\{i,j\}$-insensitive
functions.  Arguing exactly
as for Lemma~\ref{lem:isotropized} we obtain the same kind of
correlation with some s.s. extension $\tilde{\mu}$ of $\mu$ through
some factor map $\xi$, and so in light of the above nonvanishing
integral we have
\[\sfE_{\tilde{\mu}^\circ}(f_j\circ\xi -
\sfE_{\mu^\circ}(f_j\,|\,\Xi_j)\circ\xi\,|\,\tilde{\Xi}_j)\neq 0.\]
By satedness, it follows that in fact
\[\sfE_{\mu^\circ}(f_j -
\sfE_{\mu^\circ}(f_j\,|\,\Xi_j)\,|\,\Xi_j)\neq 0,\] manifestly
giving the desired contradiction. \qed

\textbf{Remark}\quad Essentially, it is the use of the
letter-replacement maps in the above proof that has been brought to
the present paper from the online project~\cite{Gow(online)}.  This
idea was brought to my attention during discussions with Terence
Tao, a more active participant in that project. \fin

We can now give our main structural result for sated laws.

\begin{thm}\label{thm:struct}
If $\mu$ is a fully sated law and $\I,\I'\subseteq \binom{[k]}{\geq
2}$ are two up-sets then the oblique $\s$-algebras $\Phi^\dag_\I$
and $\Phi^\dag_{\I'}$ are relatively independent over
$\Phi^\dag_{\I\cap\I'}$ under $\mu^\line$.
\end{thm}

\textbf{Remark}\quad This result together with
Theorem~\ref{thm:satedsexist} amount to our analog for s.s. laws of
the representation theorems for partially exchangeable
arrays~(\cite{Kal92}). \fin

We will deduce Theorem~\ref{thm:struct} result by induction using
the following special case.

\begin{lem}
If $\mu$ is a fully sated law, $\I \subseteq \binom{[k]}{\geq 2}$ is
an up-set and $e$ is a member of $\binom{[k]}{\geq 2}\setminus \I$
of maximal size then the oblique $\s$-algebras $\Phi^\dag_e$ and
$\Phi^\dag_\I$ are relatively independent over $\Phi^\dag_{\langle
e\rangle \cap\I}$ under $\mu^\line$.
\end{lem}

\textbf{Proof}\quad Suppose that $F_1 \in
L^\infty(\mu^\line|_{\Phi^\dag_e})$ and $F_2 \in
L^\infty(\mu^\line|_{\Phi^\dag_{\I}})$. It will suffice to show that
\[\int_{X^k}F_1F_2\,\d\mu^\line = \int_{X^k}\sfE_{\mu^\line}(F_1\,|\,\Phi^\dag_{\I\cap\langle
e\rangle})\cdot F_2\,\d\mu^\line.\]

Pick $i\in e$ and $f_1 \in L^\infty(\mu^\circ|_{\Phi_e})$ such that
$F_1 = f_1\circ\pi_i$ $\mu^\line$-almost surely.

Let $\{a_1,a_2,\ldots,a_q\}$ be the antichain of minimal elements in
$\I$; this clearly generates $\I$ as an up-set. Since $e \not\in\I$
we must have $a_s \setminus e\not=\emptyset$ for each $s \leq q$.
Pick $i_s \in a_s\setminus e$ arbitrarily for each $s\leq q$, so
that $\Phi^\dag_{a_s} = \pi_{i_s}^{-1}(\Phi_{a_s})$ up to
$\mu^\line$-negligible sets.

Now, since $\Phi^\dag_{\I} = \bigvee_{s\leq q}\Phi^\dag_{a_s}$,
$F_2$ may be approximated arbitrarily well in $L^1(\mu^\line)$ by
sums of products of the form $\sum_p\prod_{s\leq
q}\phi_{s,p}\circ\pi_{i_s}$ with $\phi_{s,p}\in
L^\infty(\mu^\circ|_{\Phi_{a_s}})$, and so by continuity and
linearity it suffices to assume that $F_2$ is an individual such
product term. This represents $F_2$ as a function of coordinates in
$X^k$ indexed only by members of $\{i_1,i_2,\ldots,i_q\}\subseteq
[k]\setminus e$, and now we appeal to Theorem~\ref{thm:char} to
deduce that
\begin{multline*}
\int_{X^k}F_1\cdot \prod_{s\leq q}\phi_{s,p}\circ\pi_{i_s}\,\d\mu^\line\\
=
\int_{X^k}\sfE_{\mu^\circ}\Big(f_1\,\Big|\,\bigvee_{j\in[k]\setminus
e}\Phi_{\{i,j\}}\Big)\circ\pi_i\cdot \prod_{s\leq
q}\phi_{s,p}\circ\pi_{i_s}\,\d\mu^{\line}.
\end{multline*}

However, now Lemma~\ref{lem:isotropized} and the fact that $f_1$ is
already $\Phi_e$-measurable imply that
\[\sfE_{\mu^\circ}\Big(f_1\,\Big|\,\bigvee_{j\in[k]\setminus e}\Phi_{\{i,j\}}\Big) = \sfE_{\mu^\circ}\Big(f_1\,\Big|\,\bigvee_{j\in [k]\setminus e} \Phi_{e \cup \{j\}}\Big),\]
and since each $e \cup \{j\} \in \I$ (by the maximality of $e$ in
$\cal{P}[k]\setminus\I$), under $\pi_i$ this conditional expectation
must be identified with
$\sfE_{\mu^\line}(F_1\,|\,\Phi^\dag_{\I\cap\langle e\rangle})$, as
required. \qed

\textbf{Proof of Theorem~\ref{thm:struct}}\quad We fix $\I$ and
prove this by induction on $\I'$. If $\I' \subseteq \I$ then the
result is clear, so now let $e$ be a minimal member of
$\I'\setminus\I$ of maximal size, and let $\I'' :=
\I'\setminus\{e\}$.  It will suffice to prove that if $F \in
L^\infty(\mu^\line|_{\Phi^\dag_{\I'}})$ then
\[\sfE_{\mu^\line}(F\,|\,\Phi^\dag_{\I}) = \sfE_{\mu^\line}(F\,|\,\Phi^\dag_{\I\cap\I'}),\]
and furthermore, by approximation, to do so only for $F$ that are of
the form $F_1\cdot F_2$ with $F_1 \in
L^\infty(\mu^\line|_{\Phi^\dag_{\langle e \rangle}})$ and $F_2 \in
L^\infty(\mu^\line|_{\Phi^\dag_{\I''}})$. However, for these we can
write \begin{multline*} \sfE_{\mu^\line}(F\,|\,\Phi^\dag_{\I}) =
\sfE_{\mu^\line}\big(\sfE_{\mu^\line}(F\,|\,\Phi^\dag_{\I\cup\I''})\,\big|\,\Phi^\dag_{\I}\big)\\
=
\sfE_{\mu^\line}\big(\sfE_{\mu^\line}(F_1\,|\,\Phi^\dag_{\I\cup\I''})\cdot
F_2\,\big|\,\Phi^\dag_{\I}\big),
\end{multline*}
and by the preceding lemma
\[\sfE_{\mu^\line}(F_1\,|\,\Phi^\dag_{\I\cup\I''})
= \sfE_{\mu^\line}(F_1\,|\,\Phi^\dag_{(\I\cup\I'')\cap\langle
e\rangle}).\] On the other hand $(\I\cup\I'')\cap\langle e\rangle
\subseteq \I''$ (because $\I''$ contains every subset of $[k]$ that
strictly includes $e$, since $\I'$ is an up-set), and so the
preceding lemma promises similarly that
\[\sfE_{\mu^\line}(F_1\,|\,\Phi^\dag_{(\I\cup\I'')\cap\langle
e\rangle}) = \sfE_{\mu^\line}(F_1\,|\,\Phi^\dag_{\I''}).\] Therefore
the above expression for $\sfE_{\mu^\line}(F\,|\,\Phi^\dag_{\I})$
simplifies to
\begin{multline*}
\sfE_{\mu^\line}\big(\sfE_{\mu^\line}(F_1\,|\,\Phi^\dag_{\I''})\cdot
F_2\,\big|\,\Phi^\dag_{\I}\big) =
\sfE_{\mu^\line}\big(\sfE_{\mu^\line}(F_1\cdot
F_2\,|\,\Phi^\dag_{\I''})\,\big|\,\Phi^\dag_{\I}\big)\\
=
\sfE_{\mu^\line}\big(\sfE_{\mu^\line}(F\,|\,\Phi^\dag_{\I''})\,\big|\,\Phi^\dag_{\I}\big)
= \sfE_{\mu^\line}(F\,|\,\Phi^\dag_{\I\cap \I''}) =
\sfE_{\mu^\line}(F\,|\,\Phi^\dag_{\I\cap\I'}),
\end{multline*}
by the inductive hypothesis applied to $\I''$ and $\I$, as required.
\qed

\section{The Density Hales-Jewett Theorem for sated
laws}\label{sec:final}

Since it is clear that the assertion of Theorem~\ref{thm:inf-DHJ2}
holds for a s.s. law if if holds for any extensions of that law, by
Theorem~\ref{thm:satedsexist} it suffices to prove
Theorem~\ref{thm:inf-DHJ2} in case $\mu$ is fully sated.

In this case Theorem~\ref{thm:char} and Theorem~\ref{thm:struct}
together give quite a detailed picture of the joint distribution of
the factors $\Phi^\dag_\I$ under $\mu^\line$, and it turns out that
this structure is enough to enable a proof of that theorem along the
same lines as for the multidimensional Szemer\'edi Theorem
in~\cite{Aus--newmultiSzem}. In particular,
Theorem~\ref{thm:inf-DHJ2} now follows from an `infinitary removal
lemma' essentially identical to that used
in~\cite{Aus--newmultiSzem} (Proposition 6.1 of that paper), which
was in turn based on Tao's `infinitary hypergraph removal lemma'
in~\cite{Tao07}, with some modifications to fit the context of a
proof of multiple recurrence. The version we will use below is
lifted almost verbatim from~\cite{Aus--newmultiSzem}, and is
amenable to an identical proof from Theorem~\ref{thm:struct} as for
that result from Corollary 5.2 of~\cite{Aus--newmultiSzem}, so we
only state the result here.

\begin{prop}[Infinitary removal lemma]\label{prop:removal}
Suppose that $\mu^\line$ is the line marginal of a fully sated s.s. law $\mu$,
and so has the structure described by Theorem~\ref{thm:struct}, and
that $\I_{i,j}$ for $i=1,2,\ldots,d$ and $j = 1,2,\ldots,k_i$ are
collections of up-sets in $\binom{[k]}{\geq 2}$ such that $[k] \in
\I_{i,j} \subseteq \langle i\rangle$ for each $i,j$, and suppose
further that the sets $A_{i,j}\in \Phi_{\I_{i,j}}$ are such that
\[\mu^\line\Big(\prod_{i=1}^d\Big(\bigcap_{j=1}^{k_i}A_{i,j}\Big)\Big) = 0.\]
Then we must also have
\[\mu^\circ\Big(\bigcap_{i=1}^d\bigcap_{j=1}^{k_i}A_{i,j}\Big) = 0.\]
\qed
\end{prop}

This is proved by an induction on a suitable ordering of the
possible collections of up-sets $(\I_{i,j})_{i,j}$, appealing to a
handful of different possible cases at different steps of the
induction, closely related to the induction on edge-size that
underlies the proof of the simplex removal lemma from the finitary
hypergraph regularity lemma (see, for example, Gowers~\cite{Gow06}
or Nagle, R\"odl and Schacht~\cite{NagRodSch07}).  This inductive
proof is the reason for the above statement in terms of arbitrary
collections of up-sets, but we will need only the special case
$k_i=1$, $\I_{i,1}:= \langle i\rangle$ for the proof of
Theorem~\ref{thm:inf-DHJ2}.

\textbf{Proof of Theorem~\ref{thm:inf-DHJ2} from
Proposition~\ref{prop:removal}}\quad As remarked above it suffices
to prove Theorem~\ref{thm:inf-DHJ2} for a sated law $\mu$. Given
such a law, suppose that $A \in \S$ is such that $\mu^\line(A^k) =
0$. Then by Theorem~\ref{thm:char} we have
\[\int_{X^k}\prod_{i=1}^k \sfE_{\mu^\circ}(1_A\,|\,\Phi_{\langle i\rangle})\circ\pi_i\,\d\mu^\line = \mu^\line(A^k) = 0.\]
Now the level set $B_i := \{\sfE_{\mu^\circ}(1_A\,|\,\Phi_{\langle
i\rangle})
> 0\}$ lies in $\Phi_{\langle i\rangle}$, and the above equality certainly
implies that also $\mu^\line(B_1\times B_2\times \cdots\times B_d) =
0$.  Now, on the one hand, setting $k_i=1$, $\I_{i,1}:= \langle
i\rangle$ and $A_{i,1} := B_i$ for each $i \leq d$,
Proposition~\ref{prop:removal} tells us that $\mu(B_1\cap
B_2\cap\cdots\cap B_d) = 0$, while on the other we must have
$\mu(A\setminus B_i) = 0$ for each $i$, and so overall $\mu(A) \leq
\mu(B_1\cap B_2\cap\cdots\cap B_d) + \sum_{i=1}^d\mu(A\setminus B_i)
= 0$, as required. \qed

\parskip 0pt

\bibliographystyle{abbrv}
\bibliography{newinfDHJ}

\vspace{10pt}

\small{\textsc{Department of Mathematics, Brown University,
Providence, RI 02912, USA}}

\vspace{5pt}

\small{Email: \verb|timaustin@math.brown.edu|}

\vspace{5pt}

\small{URL: \verb|http://www.math.brown.edu/~timaustin|}

\end{document}